\documentclass[11pt, oneside]{article}   	
\usepackage{amsmath,geometry}                		
\usepackage{graphicx}				
\usepackage{amssymb}
\usepackage[numbers]{natbib}
\usepackage{url}

\title{Equivalence of OEIS A007729 and A174868}
\author{Michael J. Collins\\Daniel H. Wagner Associates\\mjcollins10@gmail.com\\\\David Wilson\\davidwwilson710@gmail.com.}
\date{}		

\begin{document}
\maketitle
\begin{abstract}
We verify the conjecture that the sixth binary partition function \cite{A007729} is equal (aside from the initial zero term) to the partial sums of the Stern-Brocot sequence \cite{A174868}:  $$0,1, 2, 4, 5, 8, 10, 13, 14, 18, 21, 26, 28, 33, 36, 40, 41, 46\ldots$$.
\end{abstract}
Let $b'_k$ be the \emph{sixth binary partition function}, which is the number of ways to write $k$ as a sum
\begin{equation*}
k = \sum_{i \geq 0} \varepsilon_i 2^i
\end{equation*}
with $\varepsilon_i \in \{0,1,2,3,4,5\}$.
We obtain
\begin{equation}
b'_{2k} = b'_k + b'_{k-1} + b'_{k-2}
\end{equation}
by counting the number of representations of $2k$ with $\varepsilon_0 = 0,2$ or $4$; in each case we have a representation $\hat \varepsilon$ of  $\frac{2k-\varepsilon_0}{2}$ by taking $\hat\varepsilon_i = \varepsilon_{i+1}$, and the correspondence is clearly one-to-one. Also  $b'_{2k+1}=b'_{2k}$, since we can get a representation of $2k+1$ only by taking a representation of $2k$ and adding 1 to $\varepsilon_0$. Thus
 \begin{equation}\begin{aligned}
b'_{2k} &= 2b'_{k-1} + b'_k \ \ \ (k \mbox{ even}) \\
b'_{2k} &= 2b'_{k-1} + b'_{k-2} \ \ \ (k \mbox{ odd})
\end{aligned}\end{equation}
since $b'_{k-1} $ equals either $b'_k$ or $b'_{k-2}$. 

We eliminate the even/odd repetition by defining $b_k = b'_{2k}$. Then $b_0=1,b_1=2$ and
\begin{equation}\begin{aligned}
b_{2k} = 2b'_{2k-1} + b'_{2k} &=    2b_{k-1} + b_k\\
b_{2k+1} = 2b'_{2k} + b'_{2k-1} &= 2b_k + b_{k-1} \ .
\end{aligned}\end{equation}
This is A007729. If we prepend a zero, defining $\hat b_0=0$ and $\hat b_k = b_{k-1}$ we obtain
\begin{equation}\begin{aligned}
\hat b_{2k} &= 2 \hat b_k + \hat b_{k-1} \\
\hat b_{2k+1}  &= 2 \hat b_k + \hat b_{k+1}\ .
\end{aligned}\end{equation}
The same recurrence with the same initial conditions gives A174868, the partial sums of the Stern-Brocot sequence \cite{A002487}. The Stern-Brocot sequence itself can be defined by $s_0=0, s_1=1$, and
\begin{equation}\begin{aligned}
s_{2k} &= s_k \\
s_{2k+1} &= s_k + s_{k+1}\ .
\end{aligned}\end{equation}
The partial sums are $\sigma_k = \sum_{0 \leq i \leq k} s_k$. Letting  $\ell_j = s_{2j-1}+s_{2j} = 2s_j+s_{j-1}$ we get
\begin{equation*}
\sigma_{2k} = \sum_{1 \leq j \leq k} \ell_j  = 2\sigma_k + \sigma_{k-1}
\end{equation*}
and similarly with  $\ell'_j = s_{2j}+s_{2j+1} = 2s_j+s_{j+1}$,
\begin{equation*}
 \sigma_{2k+1} = \sum_{0 \leq j \leq k} \ell'_j  = 2\sigma_k + \sigma_{k+1} \ .
\end{equation*}

\bibliography{EquivSeq}{}
\bibliographystyle{plainnat}
\end{document}